\documentclass[12pt,twoside]{article}

\usepackage[english]{babel}
\usepackage{amsmath}
\usepackage{amsfonts,dsfont}
\usepackage{amssymb}
\usepackage{enumerate}
\usepackage{mathrsfs}
\usepackage{amssymb}
\usepackage[all]{xy}
\usepackage{graphics}
\usepackage{pstricks,pst-plot,pst-node}

\usepackage{dsfont}
\usepackage{amsthm}

\date{}

\begin{document}

\newtheorem{theorem}{\quad theorem}[section]

\newtheorem{Definition}[theorem]{\quad Definition}

\newtheorem{Corollary}[theorem]{\quad Corollary}

\newtheorem{lemma}[theorem]{\quad lemma}

\centerline{}

\centerline{}

\centerline {\Large{\bf A note on commuting graphs of matrix rings over fields.}}

\centerline{}

\centerline{\bf {C. Miguel}}

\centerline{}

\centerline{ Instituto de Telecomunuca\c
c\~oes,}

\centerline{P\'olo de Covilh\~a,}

\centerline{celino@ubi.pt}

\begin{abstract}We will give a short proof of the fact that if the algebraic  closure of a field $\mathds F$ is a finite extension, then for $n\geq 3$ the commuting graph
$\Gamma(M_n(\mathds F))$ is connected and its diameter is four.

\end{abstract}




\section{Introduction}
Let $\mathds F$ be a field and $M_n(\mathds F)$ be the algebra of all $n\times n$ matrices over $\mathds F$. The commuting graph of $M_n(\mathds F)$, denoted by $\Gamma(M_n(\mathds F))$,
is the graph whose vertices are all non-central matrices and two distinct vertices $A$ and $B$ are adjacent if and only if $AB=BA$.

The connectedness and diameter of the commuting  $\Gamma(M_n(\mathds F))$  have
been  studied extensively; see \cite{ak2,dol, do,gi}. Note that for $n=2$ the commuting graph $\Gamma(M_n(\mathds F))$ is always disconnected; see \cite[remark 8]{ak4}. So, from now on we assume that $n$ is an integer greater than $2$.   If the field $\mathds F$ is algebraically closed, then the commuting graph $\Gamma(M_n(\mathds F))$ is connected and its diameter is always equal to four; see \cite{ak2}. If the field $\mathds F$ is not algebraically closed, then the commuting graph $\Gamma(M_n(\mathds F))$ may be disconnected for an arbitrarily large integer $n$; see \cite{ak3}. However, in $\cite{ak2}$ S. Akbari et al.  proved that for any field $\mathds F$ if the graph $\Gamma(M_n(\mathds F))$ is connected, then its diameter is between four and six. Also in $\cite{ak2}$ the authors conjectured that the upper bounded six can be improved to five.
Recently, Yaroslav Shitov in \cite{sh} gave an example of a field $\mathds F$ such that $\Gamma(M_{38}(\mathds F))$ is connected with diameter six, disproving the conjecture of S. Akbari et al.

The aim of this note is to give a short proof that if the algebraic closure $\overline{\mathds F}$ of the field $\mathds F$ is a finite extension of $\mathds F$, then the commuting graph $\Gamma(M_n(\mathds F))$ is connected and its diameter is four.
We will prove the following:

\begin{theorem}\label{tf} If the algebraic closure $\overline{\mathds F}$ of a field $\mathds F$ is a finite extension of $\mathds F$,
then
the commuting graph $\Gamma(M_n(\mathds F))$ is connected and its diameter is four. \end{theorem}

Notice that, the well-known fact that for matrix rings over algebraically closed fields, if the commuting graph is connected, then its diameter ir four, is a particulary case of theorem \ref{tf}. Moreover, the results proved in \cite{mi} and \cite{ga} for the field $\mathds R$ of real numbers can also be deduced from
theorem \ref{tf}.

\section{Preliminaries.}

Our proof involves looking at the idempotents and nilpotent matrices and apply results from the structure theory of rings.
Recall that in a ring $R$ both the zero and the identity are always idempotent. An indempotent wich is different from these is called non trivial idempotent. Note that, For a field
$\mathds F$ the central idempotent in the matrix ring $M_n(\mathds F)$ are exactly the trivial idempotent.

\begin{lemma}\label{lf} Let $\mathds F$ be field and assume that all the irreducible polynomials in $\mathds F[x]$ have
degree $< n$. Then, every matrix $A\in M_n(\mathds F)$ commute with a non zero nilpotent matrix or with a nontrivial idempotent matrix.\end{lemma}

{\bf Proof.} First observe that, if the matrix $A$ is  derrogatory, then by the rational canonical form we can find  a nonsingular matrix $S\in M_n(\mathds F)$ such that
\begin{equation*}A=S^{-1}(C_1\oplus\ldots\oplus C_k)S,\end{equation*}
where $k\geq 2$ and $C_i$, for $i=1,\ldots ,k$, is a companion matrix of a monic polynomial over $\mathds F$. Hence, $A$ commutes with the non trivial idempotent
\begin{equation*}E=S^{-1}(I_{C_1}\oplus 0_{C_2}\oplus\ldots\oplus 0_{C_k})S,\end{equation*}
where $I_{C_1}$ is the identity matrix whose order equals the order of $C_1$ and $0_{C_i}$, for $i=2, \ldots, k$, is the zero matrix whose order equals the order of $C_i$.

Now, assume that the matrix $A$ is non derrogatory and denote by $\langle A\rangle$ the $\mathds F$-subalgebra of $M_n(\mathds F)$ generated by $A$. Since the $\mathds F$-subalgebra
$\langle A\rangle$ is artinian it follows that the Jacobson radical $\mathcal{J}(\langle A\rangle)$ is nilpotent; see \cite[p.658]{lang}. Consequently, if $\langle A\rangle$ is not nilsemisimple it necessarily
contains a non zero nilpotent matrix $N$ and we get the result.

We are now reduced to the case where the matrix $A$ is non derrogatory and the $\mathds F$-algebra $\langle A\rangle$ is semisimple.
Since a semisimple artinian commutative ring is isomorphic to sum of fields; see \cite[p.661]{lang}, if follows that
\begin{equation}\label{eq1}\langle A\rangle\cong\mathds F_1\oplus\ldots\oplus\mathds F_s,\end{equation}
where each $\mathds F_i$, for $i=1\ldots, s$, is a field.
If in equation (\ref{eq1}) we have $s\geq 2$, then  $\mathds F_1\oplus\ldots\oplus\mathds F_s$ contains the non-trivial   idempotent $(1, 0, \ldots , 0)$. Let $\psi:F_1\oplus\ldots\oplus\mathds F_s\rightarrow \langle A\rangle $ be an isomorphism. Clearly, a non-trivial idempotent of $\mathds F_1\oplus\ldots\oplus\mathds F_s$ is mapped
by $\psi$ to a non-trivial idempotent in $\langle A\rangle$.
 As we have noted before, in the matrix ring $M_n(\mathds F)$ the central idempotent are exactly the trivial idempotent. Hence, the idempotent $\psi(1, 0, \ldots , 0)$ is not central and we get the result.

 We conclude the proof by showing that in equation (\ref{eq1}) we necessary have $s\geq 2$. Indeed, if we had $s=1$ in equation (\ref{eq1}), then
 the subalgebra  $\langle A\rangle$ is
itself a field and by the first isomorphism theorem for rings we obtain
\begin{equation*}\langle A\rangle\cong\mathds F[x]/( p_A),\end{equation*}
where $( p_A)$ is the ideal, in the polynomial ring $\mathds F[x]$, generated by the minimum polynomial $p_A$ of $A$, over the field $\mathds F$. Since the ideal $( p_A)$ is maximal it follows that $p_A$ is irreducible over $\mathds F$. Now the fact that $A$ is non derrogatory implies that $p_A$ has degree $n$, contrary to hypothesis. $\Box$

\

 Observe that if a matrix $A$ commutes with  a noncentral nilpotent matrix $N$ with index of  nilpotency  $\alpha$, that is, $\alpha=min\{i\in\mathds N : N^i=0\}$, then the matrix $A$ also commutes with $N^{\alpha-1}$. Consequently, if a matrix commutes with a noncentral nilpotent matrix $N$ we may assume that its index of nilpotency is $2$.

\section{Proof of the main theorem.}

We can now prove our main result.

\

{\bf Proof of Theorem \ref{tf}.} Following \cite[theorem\;11]{ak2}, if $E_1, E_2\in M_n(\mathds F)$ are two noncentral idempotent matrices, then $d(E_1, E_2)\leq 2$ in the commuting graph $\Gamma(M_n(\mathds F))$. Also by  \cite[theorem\;9]{ak2}
 if $N_1, N_2\in M_n(\mathds F)$ are two non-zero nilpotent matrices, both of index of nilpotence $2$, then $d(N_1, N_2)\leq 2$ in $\Gamma(M_n(\mathds F))$. Finally, from \cite[proposition 1]{ga} if $E, N\in M_n(\mathds F)$ are  such that $E$ is idempotent and $N$ is nilpotent of index of nilpotence $2$, then $d(E, N)\leq 2$ in $\Gamma(M_n(\mathds F))$.
  By  combining these three results with lemma \ref{lf} we get the result. $\Box$

\

{\bf Final Remarks:} For a field $\mathds F$ if the commuting graph $\Gamma(M_n(\mathds F))$ is connected,  then there are three possibilities
for the diameter namely $4$, $5$ or $6$. According to theorem \ref{tf}  for fields $\mathds F$ with an algebraic closure that is a finite extension the graph
$\Gamma(M_n(\mathds F))$ is connected with diameter $4$.  On the other hand, Yaroslav Shitov in \cite{sh} gave an example of a field $\mathds F$ such that $\Gamma(M_{38}(\mathds F))$ is connected with diameter $6$. However, there is no known example of a connected commuting graph of a matrix ring over a field with diameter $5$. The question natural arises
  whether a field $\mathds F$ exists such that for some integer $n$ the commuting $\Gamma(M_n(\mathds F))$ is connected with diameter five?

\

{\bf Acknowledgments.}

This work was supported by FCT project UID/EEA/50008/2013.

\end{document}